\documentclass[reqno,12pt]{amsart}

\newtheorem{theorem}{Theorem}[section]
\newtheorem{lemma}[theorem]{Lemma}

\newtheorem{proposition}[theorem]{Proposition}

\theoremstyle{definition}
\newtheorem{definition}[theorem]{Definition}

\theoremstyle{remark}
\newtheorem{remark}[theorem]{Remark}

\newcommand{\mysection}[1]{\section{#1}
\setcounter{equation}{0}}

\newcommand{\bR}{\mathbb R}

\renewcommand{\epsilon}{\varepsilon}

\begin{document}
\title[critical quasi-geostrophic equations
]{Global well-posedness and a decay estimate for the critical 
dissipative quasi-geostrophic equation in the whole space}

\author[H. Dong]{Hongjie Dong}
\address[H. Dong]
{School of Mathematics, Institute for Advanced Study,
Einstein Drive, Princeton, NJ 08540, USA}
\email{hjdong@math.ias.edu}

\author[D. Du]{Dapeng Du}
\address[D. Du]{School of Mathematical Sciences,
Fudan University, Shanghai 200433, P.R. China}
\email{dpdu@fudan.edu.cn}

\thanks{H. D. is partially supported by the National Science
Foundation under agreement No. DMS-0111298. D. D. is partially
supported by a postdoctoral grant from School of Mathematical
Sciences at Fudan University.}



\begin{abstract}
We study the critical dissipative quasi-geostrophic equations in
$\bR^2$ with arbitrary $H^1$ initial data. After showing certain
decay estimate, a global well-posedness result is proved by adapting
the method in \cite{kiselev} with a suitable modification. A decay
in time estimate for higher order homogeneous Sobolev norms of
solutions is also discussed.
\end{abstract}

\maketitle

\mysection{Introduction}
                                                    \label{intro}

In this note, we consider the initial value problem of 2D
dissipative quasi-geostrophic equations
\begin{equation}
                                        \label{qgeq1}
\left\{\begin{array}{l l}
\theta_t+u\cdot \nabla\theta+(-\Delta)^{\gamma/2}\theta=0
\quad &\text{on}\, \bR^2\times (0,\infty),\\
\theta(0,x)=\theta_0(x)\quad & x\in \bR^2,\end{array}\right.
\end{equation}
where $\gamma\in (0,2]$ is a fixed parameter and the velocity
$u=(u_1,u_2)$ is divergence free and determined by the Riesz
transforms of the potential temperature $\theta$:
$$
u=(-{\mathcal R}_2\theta,{\mathcal R}_1\theta)=(-\partial_{x_2}
(-\Delta)^{1/2}\theta, \partial_{x_1}(-\Delta)^{1/2}\theta).
$$
The main problem addressed here is the global regularity of
\eqref{qgeq1} with $\gamma=1$ and arbitrary $H^1$ initial data.

Equation \eqref{qgeq1} is an important model in geophysical fluid
dynamics. It is derived from general quasi-geostrophic equations in
the special case of constant potential vorticity and buoyancy
frequency. When $\gamma=1$, it is the dimensionally correct analogue
of the 3D incompressible Navier-Stokes equations. The critical
dissipative quasi-geostrophic equations is an interesting model for
investigating existence issues on genuine 3D Navier-Stokes
equations. Recently, this equation has been studied intensively, see
\cite{cordoba}, \cite{const1}, \cite{const2}, \cite{const3},
\cite{ju1}, \cite{ju2}, \cite{miao}, \cite{miura}, \cite{resnick},
\cite{wu1}, \cite{wu2}, \cite{wu3}, \cite{yu} and references
therein.

The cases $\gamma>1, \gamma=1$ and $\gamma<1$ are  called
sub-critical, critical and super-critical respectively. The
sub-critical case is better understood and the global well-posedness
result is well-known. For this case, we refer the readers to Wu
\cite{wu1}, Carrilo and Ferreira \cite{Carrilo}, Constantin and Wu
\cite{const2}, Dong and Li \cite{dongli} and reference therein.

The cases of critical and super-critical dissipative
quasi-geostrophic equations still have quite a few unsolved
problems. One major problem is the issue of global regularity or
breakdown of regular solutions, which was suggested by Klainerman
\cite{Klai} seven years ago as one of the most challenging PDE
problems of the twenty-first Century. In the critical case,
Constantin, C\'ordoba and Wu \cite{const1} gave a construction of
global regular solutions for initial data in $H^1$ under a smallness
assumption of $L^\infty$ norm of the data. For other results about
local well-posedness and small-data global well-posedness  in
various function spaces, see also Chae and Lee \cite{chae}, Ju
\cite{ju1}, \cite{ju3}, Miura \cite{miura} and reference therein.
For the issue of global regularity with large initial data,
breakthrough only occurred recently. Caffarelli and Vasseur
\cite{caf} constructed a global regular Leray-Hopf type weak
solution for the critical quasi-geostrophic equations with merely
$L^2$ initial data. The global well-posedness for the critical
quasi-geostrophic equations with periodic $C^\infty$ data was proved
by Kiselev, Nazarov and Volberg in an elegant paper \cite{kiselev}.
Their argument is based on a certain non-local maximum principle for
a suitable chosen modulus of continuity.

Miura \cite{miura} recently established the local in time existence
of a unique regular solution for large initial data in the critical
Sobolev space $H^{2-\gamma}$. A similar result was also obtain
independently in Ju \cite{ju3} by using a different approach. Very
recently, the first author showed that the solutions by Miura and Ju
have higher regularity and are global in time with periodic $H^1$
data. Roughly speaking, it is proved that the smoothing effect of
the equations in spaces is the same for the corresponding linear
equations. For other results about the critical and super-critical
dissipative quasi-geostrophic equations, we also refer the readers
to \cite{const3}, \cite{ju2}, \cite{miao}, \cite{wu2}, \cite{wu3}
and \cite{yu}. However, at present the following aforementioned
problem suggested by Klainerman is still open:

\noindent {\em For $\gamma=1$, is \eqref{qgeq1} globally well-posed
with arbitrary smooth in $\bR^2$ ?}

We give an affirmative answer to this question (Theorem \ref{thm1})
with initial data in the critical Sobolev space $H^1$. Our strategy
is to apply a local smoothing result proved in \cite{dong} and then
adapt the idea in \cite{kiselev} with a proper modification. One
essential difference between periodic and non-periodic settings is
that in periodic domains, one can appeal certain compactness
property, which is not valid in the whole space. We circumvent this
difficult by showing some decay estimate of solutions as space
variables go to infinity. Moreover, we show that the solution is
actually a smooth classical solution to \eqref{qgeq1} and higher
order homogeneous Sobolev norms of the solution decay polynomially
as $t$ goes to infinity (Theorem \ref{thm2}).

The remaining part of the note is organized as follows: after
reviewing some local well-posedness and smoothing results, the main
theorem is given in the following section. We prove the theorem in
Section \ref{proofsec}. Section \ref{remarksec} is devoted to a
decay in time estimate of higher order homogeneous Sobolev norms of
solutions.

\mysection{Main Theorem}

The local well-posedness of \eqref{qgeq1} with $H^{2-\gamma}$ data
is recently established by Miura and Ju independently.

\begin{proposition}
            \label{prop0}
Let $\gamma\in (0,1]$ and $\theta_0\in H^{2-\gamma}(\bR^2)$. Then
there exists $T>0$ such that the initial value problem for
\eqref{qgeq1} has a unique solution
\begin{equation*}
\theta(t,x)\in C([0,T); H^{2-\gamma}(\bR^2))\cap L^2(0,T;
H^{2-\gamma/2}(\bR^2)).
\end{equation*}
Moreover the solution $\theta$ satisfies
\begin{equation*}
\sup_{0<t<T}t^{\beta/\gamma}\|\theta(t,\cdot)\|_{\dot
H^{2-\gamma+\beta}}<\infty,
\end{equation*}
for any $\beta\in [0,\gamma)$ and
\begin{equation*}
\lim_{t\to 0}t^{\beta/\gamma}\|\theta(t,\cdot)\|_{\dot
H^{2-\gamma+\beta}}=0,
\end{equation*}
for any $\beta\in (0,\gamma)$.
\end{proposition}

The following proposition is the main result of \cite{dong}, which
says that the solution in Proposition \ref{prop0} has higher
regularities.

\begin{proposition}
                    \label{prop1}
The solution $\theta$ in Proposition \ref{prop0} satisfies
\begin{equation}
            \label{eq2.08b}
\sup_{0<t<T}t^{\beta/\gamma}\|\theta(t,\cdot)\|_{\dot
H^{2-\gamma+\beta}}<\infty,
\end{equation}
for any $\beta\geq 0$ and
\begin{equation}
            \label{eq2.10b}
\lim_{t\to 0}t^{\beta/\gamma}\|\theta(t,\cdot)\|_{\dot
H^{2-\gamma+\beta}}=0,
\end{equation}
for any $\beta>0$.
\end{proposition}

By the Sobolev imbedding theorem, the previous proposition implies
that the solution in Proposition \ref{prop0} is infinitely
differentiable in $x$ with bounded derivatives for any $t\in (0,T)$.
Then because of the first equation in \eqref{qgeq1}, it is
infinitely differentiable in both $x$ and $t$ with bounded
derivatives for any $t\in (0,T)$. Therefore, $\theta$ is actually a
classical solution of \eqref{qgeq1}.

In the sequel, we always assume $\gamma=1$, i.e. the critical case.
Next we state our main theorem.

\begin{theorem}
                                \label{thm1}
With $\theta_0\in H^1(\bR^2)$, the initial value problem
\eqref{qgeq1} has a unique global solution
\begin{equation}
                                            \label{eq4.18}
\theta(t,x)\in C_b([0,\infty); H^1(\bR^2))\cap L^2(0,\infty;
H^{3/2}(\bR^2)).
\end{equation}
\end{theorem}

\begin{remark}
                                \label{rem1}
In \cite{dong}, it is shown that with zero-mean periodic $H^1$ data,
the solution and all its derivatives decay exponentially as $t$ goes
to infinity. This is certainly not the case for equations with
non-periodic data. Moreover, in the periodic setting
$\theta(t,\cdot)$ is spatial periodic for sufficiently large $t$. We
conjecture this is still true for the solution in Theorem
\ref{thm1}.
\end{remark}


\mysection{Proof of Theorem \ref{thm1}}
                \label{proofsec}

This section is devoted to the proof of Theorem \ref{thm1}. The
argument is mainly based on a decay estimate of $|\nabla_x \theta|$
as $x\to \infty$ and a non-local maximum principle as in
\cite{kiselev} with a proper modification.

\begin{definition}
We say a function $f$ has modulus of continuity $\omega$ if
$|f(x)-f(y)|\leq \omega(|x-y|)$, where $\omega$ is an unbounded
increasing continuous concave function $\omega$ : $[0,+\infty)\to
[0,+\infty)$. We say $f$ has strict modulus of continuity $\omega$
if the inequality is strict for $x\neq y$.
\end{definition}

Recall the remark before Theorem \ref{thm1}. After fixing a time
$t_1\in (0,T)$ and considering $\theta(t-t_1)$ instead of $\theta$,
we may assume $\theta_0\in H^1\cap C^\infty$ and $\theta_0$ is
bounded along with all its derivatives. Therefore, for any unbounded
continuous concave function $\omega$ satisfying
\begin{equation}
                    \label{eq10.18}
\omega(0)=0,\quad \omega'>0,\quad \omega'(0)<+\infty,\quad
\lim_{\xi\to 0^+} \omega''(\xi)=-\infty,
\end{equation}
we can
find a constant $C>0$ such that $\omega(\xi)$ is a {\em strict}
modulus of continuity of $\theta_0(Cx)$. Due
to the scaling property of \eqref{qgeq1},
$\theta_c(t,x)=\theta(Ct,Cx)$ is the solution of \eqref{qgeq1} with
initial data $\theta_0(Cx)$. Thus if we can show that suitable
modulus of continuity $\omega$ is preserved by the dissipative
evolution so that $\theta_c$ is a global solution, the same is true
for $\theta$.

With aforementioned property \eqref{eq10.18} of $w$, if $f\in C^2$
has modulus of continuity $\omega$, it is easy to show that
pointwisely $|\nabla f(x)|<\omega'(0)$ (see \cite{kiselev}). Assume
further that $f\in H^4(\bR^2)$, due to the Sobolev imbedding
theorem, $\nabla f(x)$ is a uniformly continuous function and goes
to zero as $|x|\to \infty$. Thus, we get
\begin{lemma}
                \label{lem3.1}
If $f\in  H^4(\bR^2)$ has modulus of continuity $\omega$ satisfying
\eqref{eq10.18}, we have $\|\nabla f\|_{L^\infty}<\omega'(0)$.
\end{lemma}

Next we show that strict modulus of continuity is preserved at least
for a short time.

\begin{lemma}
                \label{lem3.2}
Assume $\theta(t,\cdot)$ has strict modulus of continuity $\omega$
for all $t\in [0,T_1]$. Then there exists $\delta>0$ such that
$\theta(t,\cdot)$ has strict modulus of continuity $\omega$ for all
$t\in [0,T_1+\delta)$.
\end{lemma}
\begin{proof}
By the assumption,
$$
\sup_{0\leq t\leq T_1}\|\theta(t,\cdot)\|_{H^k_x}<\infty
$$
for any $k\geq 0$ and $\theta$ is smooth up to time $T_1$. Owing to
the local existence and regularity theorem, there exists a number
$\delta_1\in (0,1)$, such that we can continue $\theta$ up to time
$T_1+\delta_1$ and
$$
\sup_{0\leq t\leq
T_1+\delta_1}\|\theta(t,\cdot)\|_{H^{20}_x}<\infty.
$$
By the Sobolev imbedding
theorem and the first equation in \eqref{qgeq1}, $\theta\in
C^2([0,T_1+\delta_1]\times \bR^2)$ with bounded derivatives up to
order two. Since $\omega$ is unbounded and by Lemma \ref{lem3.1},
there exists $\delta_2\in (0,\delta_1)$ so that
\begin{equation}
                    \label{eq11.21}
|\theta(t,x)-\theta(t,y)|<\omega(|x-y|)
\end{equation}
for any $t\in [T_1,T_1+\delta_2]$ and $|x-y|\in (0,\delta_2]
\cup [\delta_2^{-1},\infty)$.

In what follows we always assume  $|x-y|\in
(\delta_2,\delta_2^{-1})$. Note that in $[T_1,T_1+\delta_2]\times
\bR^2$ $|\nabla_x \theta|$ is uniformly continuous and belongs to
$L^2([T_1,T_1+\delta_2]\times \bR^2)$. Thus, it goes to zero as
$|x|\to \infty$ uniformly in $t$, and we can find a constant  $N>0$
such that
$$
|\nabla_x \theta(t,x)|<\delta_2\omega(\delta_2^{-1})
$$
for any $t\in [T_1,T_1+\delta_1]$ and any $x$ satisfying $|x|\geq
N$. Now by the concavity of $\omega$, it is easily seen that
\eqref{eq11.21} holds for any $x,y$ satisfying $|x|\geq N$, $|y|\geq
N$. Finally, if either $|x|<N$ or $|y|<N$, we must have $ (x,y)\in
\Omega$, where $$\Omega:=\{(x,y)\in \bR^2\times
\bR^2\,|\,\max[|x|,|y|] \leq N+\delta_2^{-1},|x-y|\geq
\delta_2^{-1}\}.
$$
Because of this and since $\theta\in (T_1,\cdot)$ has strict
modulus of continuity $\omega$ and $\theta$ is a $C^2$ function,
there exists $\delta\in (0,\delta_2)$ so that \eqref{eq11.21} in
$[T_1,T_1+\delta)\times \Omega$. The lemma is proved.
\end{proof}

Notice that if $\theta(t,\cdot)$ has strict modulus of continuity
$\omega$ for all $t\in [0,T_1)$, then $\theta$ is smooth up to $T_1$
and $\theta(T_1,\cdot)$ has modulus of continuity $\omega$ by
continuity. Therefore, to show that the modulus of continuity is
preserved for all the time, it suffices to rule out the case that
\begin{equation*}
\sup_{x\neq y}\frac {\theta(T_1,x)-\theta(T_1,y)}{\omega(|x-y|)}=1.
\end{equation*}
\begin{lemma}
                    \label{lem3.4}
Under the conditions above, there exist two different points $x,y\in
\bR^2$ satisfying
$$
\theta(T_1,x)-\theta(T_1,y)=\omega(|x-y|).
$$
\end{lemma}
Assume for a moment that Lemma \ref{lem3.4} is proved. We choose a
suitable $\omega$ by letting for $r>0$
$$
\omega''(r)=-\frac{\delta_3}{r^{1/2}+r^2\log r},\quad
\omega'(r)=-\int_r^\infty \omega''(s)\,ds,\quad \omega(0)=0.
$$
Then following the argument in \cite{kiselev}, for a sufficiently
small $\delta_3$, we reach a contradiction:
$$
\frac {\partial}{\partial t}(\theta(T_1,x)-\theta(T_2,y))<0.
$$
Instead of rewriting the proof, we refer the readers to
\cite{kiselev}, where a slightly different modulus of continuity is
constructed.

Now it remains to prove Lemma \ref{lem3.4}.
\begin{proof}
By the same argument as in the proof of Lemma \ref{lem3.2}, there
exist $\delta_2,N>0$ such that \eqref{eq11.21} holds in
$\bR^2\setminus \Omega$. Then the lemma follows from the compactness
of $\Omega$.
\end{proof}

We have shown that $\omega$ is a strict  modulus of continuity for
all the time and $\theta$ doesn't have gradient blow-up.
Consequently, $\theta\in L^\infty_{\text{loc}}((0,T),H^k(\bR^2))$
for any $k\geq 0$ and it is a smooth global solution to
\eqref{qgeq1}. Due to the boundedness of the Riesz transforms and
the Sobolev imbedding theorem, $u$ is also smooth. The uniqueness
then follows in a standard way from the local uniqueness result
(see, e.g. \cite{miura}).

In the last part of this section, we shall prove \eqref{eq4.18}. By
using Theorem 4.1 of C\'ordoba and C\'ordoba \cite{cordoba} we
obtain the following decay estimate for $L^\infty$ norm of the
solution.
\begin{lemma}
                                                    \label{lem11.26}
Under the assumptions of Theorem \ref{thm1}, there exists a
positive constant $C$ depending only on $\theta_0$ so that
$$
\|\theta(t,\cdot)\|_{L^\infty}\leq \frac C {1+t}
$$ for any $t\geq t_1$.
\end{lemma}

After multiplying the first equation of \eqref{qgeq1} by  $\Delta
\theta$, integrating by parts and using the boundedness of Riesz
transforms and the Gagliardo-Nirenberg inequality, we get (see
Theorem 2.5 in \cite{const1} for details)
\begin{equation}
                \label{eq2.42}
\frac{1}{2}\frac{d}{dt}\|\nabla\theta(t,\cdot)\|_{L^2_x}^2+
\|\Lambda^{3/2}\theta(t,\cdot)\|_{L_x^2}^2
\leq C_1\|\theta(t,\cdot)\|_{L_x^\infty}\|\Lambda^{3/2}
\theta(t,\cdot)\|_{L_x^2}^2,
\end{equation}
for some constant $C_1>0$. By Lemma \ref{lem11.26}, there  exists
$T_2\geq t_1$ such that
$$
\|\theta(t,\cdot)\|_{L^\infty}\leq \frac 1 {2C_1}
$$ holds for any $t\geq T_2$.
The inequality above and \eqref{eq2.42} yield
$\|\nabla\theta(t,\cdot)\|_{L^2_x}$ is non-increasing in
$[T_2,\infty)$ and $\|\Lambda^{3/2}\theta(t,\cdot)\|_{L_x^2}\in
L^2([T_2,\infty))$. This together with the $L^p$-maximum principle
for the quasi-geostrophic equations completes the proof of Theorem
\ref{thm1}.

\mysection{A decay estimate}
                    \label{remarksec}

The last section is for a decay estimate of the solution. We show
that higher order homogeneous Sobolev norms of $\theta$ decay in
time polynomially. A similar decay estimate under a smallness
assumption can be found in \cite{dong}.

By multiplying the first equation of \eqref{qgeq1} by $\Delta^2
\theta$ instead of $\Delta \theta$, in a same fashion as
\eqref{eq2.42} one can get
\begin{equation}
                \label{eq2.43}
\frac{1}{2}\frac{d}{dt}\|\Delta\theta(t,\cdot)\|_{L^2_x}^2
+\|\Lambda^{5/2}\theta(t,\cdot)\|_{L_x^2}^2\leq
C_2\|\theta(t,\cdot)\|_{L_x^\infty}\|\Lambda^{5/2}\theta(t,\cdot)\|_{L_x^2}^2,
\end{equation}
for some constant $C_2>0$. By Lemma \ref{lem11.26},  there exists
$T_3\geq t_1$ such that
$$
\|\theta(t,\cdot)\|_{L^\infty}\leq \frac 1 {2C_2}
$$ holds for any $t\geq T_3$.
Then \eqref{eq2.43} implies that $\|\Delta\theta(t,\cdot)\|_{L^2_x}$
is non-increasing in $[T_3,\infty)$ and thus bounded in
$[1,\infty)$. Consequently, by interpolation we have
\begin{equation}
            \label{eq12.09}
\|\theta(t,\cdot)\|_{\dot H^{1+\epsilon}_x}\in L^\infty([1,\infty))
\end{equation}
for any $\epsilon\in [0,1]$.  Then we can adapt the method in
\cite{dong} to obtain the following decay in time estimate.
\begin{theorem}
                                \label{thm2}
For any $\epsilon\in (0,1)$ and any $\beta\geq \epsilon$, we have
\begin{equation}
            \label{eq10.07}
\sup_{1\leq t<\infty}t^{\beta-\epsilon}\|\theta(t,\cdot)\|_{\dot
H^{1+\beta}}<\infty,
\end{equation}
\end{theorem}

Indeed by using a commutator estimate in Lemma 2.5 of \cite{dong},
without much more work one can bootstrap from the boundedness
\eqref{eq12.09} to get \eqref{eq10.07}. Since the proof essentially
follows that of Theorem 1.3 in \cite{dong}, we omit the detail and
leave it to interested readers.

\end{document}